\input ./style/arxiv-vmsta.cfg
\documentclass[numbers,compress,v1.0.1]{vmsta}

\volume{3}
\issue{1}
\pubyear{2016}
\firstpage{47}
\lastpage{57}
\doi{10.15559/16-VMSTA50}

\usepackage{cleveref}
\DeclareMathOperator{\M}{\mathbf E}
\DeclareMathOperator{\cov}{\mathbf{cov}}
\DeclareMathOperator{\tr}{\mathsf T}
\DeclareMathOperator{\I}{\mathrm I}
\newcommand{\crefrangeconjunction}{--}
\newcommand{\prob}{\stackrel{\text{\rm P}}{\longrightarrow}}
\renewcommand{\theequation}{\arabic{section}.\arabic{equation}}
\numberwithin{equation}{section}
\crefname{equation}{equation}{}
\newtheorem{thm}{Theorem}
\newtheorem{lemma}[thm]{Lemma}
\newtheorem{cor}[thm]{Corollary}

\theoremstyle{definition}
\newtheorem{defin}[thm]{Definition}
\newtheorem{rem}[thm]{Remark}
\renewcommand{\theenumi}{\roman{enumi}}
\renewcommand{\labelenumi}{(\roman{enumi})}

\startlocaldefs
\newcommand{\lleft}{\left}
\newcommand{\rright}{\right}

\newcommand{\mathbi}[1]{\mbox{\boldmath$#1$}}
\allowdisplaybreaks
\endlocaldefs

\begin{document}
\begin{frontmatter}

\title{Asymptotic normality of total least squares estimator in a multivariate errors-in-variables model \mathbi{AX = B}}

\author{\inits{A.}\fnm{Alexander}\snm{Kukush}\corref{cor1}}\email{alexander\_
kukush@univ.kiev.ua}
\cortext[cor1]{Corresponding author.}
\author{\inits{Ya.}\fnm{Yaroslav}\snm{Tsaregorodtsev}}\email{777Tsar777@mail.ru}
\address{Taras Shevchenko National University of Kyiv, Kyiv, Ukraine}

\markboth{A. Kukush, Ya. Tsaregorodtsev}{Asymptotic normality of total least squares estimator}

\begin{abstract}
We consider a multivariate functional measurement error model
\hbox{$AX\approx B$}. The errors in $[A,B]$ are uncorrelated, row-wise
independent, and have equal (unknown) variances. We study the total
least squares estimator of $X$, which, in the case of normal errors,
coincides with the maximum likelihood one. We give conditions for
asymptotic normality of the estimator when the number of rows in $A$ is
increasing. Under mild assumptions, the covariance structure of the
limit Gaussian random matrix is nonsingular. For normal errors, the
results can be used to construct an~asymptotic confidence interval for
a linear functional of $X$.
\end{abstract}

\received{11 February 2016}
%
\revised{7 March 2016}
\accepted{11 March 2016}
\publishedonline{29 March 2016}
\end{frontmatter}

\section{Introduction}\label{}
We deal with overdetermined system of linear equations $AX\approx B$,
which is common in linear parameter estimation problem \cite{vaw}. If
the data matrix $A$ and observation matrix $B$ are contaminated with
errors, and all the errors are uncorrelated and have equal variances,
then the total least squares (TLS) technique is appropriate for solving
this system \cite{vaw}. Kukush and Van Huffel \cite{kukva} showed the
statistical consistency of the TLS estimator $\hat{X}_{\mathit{tls}}$ as the
number~$m$ of rows in $A$ grows, provided that the errors in $[A,B]$
are row-wise i.i.d. with zero mean and covariance matrix proportional
to a unit matrix; the covariance matrix was assumed to be known up to a
factor of proportionality; the true input matrix $A_0$ was supposed to
be nonrandom. In fact, in \cite{kukva} a more general, element-wise
weighted TLS estimator was studied, where the errors in $[A,B]$ were
row-wise independent, but within each row, the entries could be
observed without errors, and, additionally, the error covariance matrix
could differ from row to row. In \cite{ma}, an iterative numerical
procedure was developed to compute the elementwise-weighted TLS
estimator, and the rate of convergence of the procedure was established.

In a univariate case where $B$ and $X$ are column vectors, the
asymptotic normality of $\hat{X}_{\mathit{tls}}$ was shown by Gallo \cite{gal}
as $m$ grows. In \cite{pes}, that result was extended to mixing error
sequences. Both \cite{gal} and \cite{pes} utilized an explicit form of
the TLS solution.

In the present paper, we extend the Gallo's asymptotic normality result
to a multivariate case, where $A$, $X$, and $B$ are matrices.

Now a closed-form solution is unavailable, and we work instead with the
cost function. More precisely, we deal with the estimating function,
which is a matrix derivative of the cost function. In fact, we show
that under mild conditions, the normalized estimator converges in
distribution to a Gaussian random matrix with nonsingular covariance
structure. For normal errors, the latter structure can be estimated
consistently based on the observed matrix $[A,B]$. The results can be
used to construct the asymptotic confidence ellipsoid for a vector
$Xu$, where $u$ is a column vector of the corresponding dimension.

The paper is organized as follows. In Section~\ref{s:2}, we describe
the model, refer to the consistency result for the estimator, and
present the objective function and corresponding matrix estimating
function. In Section~\ref{s:3}, we state the asymptotic normality of
$\hat{X}_{\mathit{tls}}$ and provide a nonsingular covariance structure for a
limit random matrix. The latter structure depends continuously on some
nuisance parameters of the model, and we derive consistent estimators
for those parameters. Section~\ref{s:4} concludes. The proofs are given
in Appendix. There we work with the estimating function and derive an
expansion for the normalized estimator using Taylor's formula. The
expansion holds with probability tending to $1$.

Throughout the paper, all vectors are column ones, $\M$ stands for the
expectation and acts as an operator on the total product, $\cov(x)$
denotes the covariance matrix of a random vector $x$, and for a
sequence of random matrices $\{X_m,m\geq1\}$ of the same size, the
notation $X_m=O_p(1)$ means that the sequence $\{\|X_m\|\}$ is
stochastically bounded, and $X_m=o_p(1)$ means that $\|X_m\|\prob0$.
By $\I_p$ we denote the unit matrix of size $p$.
\section{Model, objective, and estimating}\label{s:2}
\subsection{The TLS problem}
Consider the model $AX\approx B$. Here $A\in\mathbb{R}^{m\times n}$ and
$B\in\mathbb{R}^{m\times d}$ are observations, and $X\in\mathbb
{R}^{n\times d}$ is a parameter of interest. Assume that
\begin{equation}
\label{2.1} A=A_0+\tilde{A}, \qquad B=B_0+\tilde{B},
\end{equation}
and that there exists $X_0\in\mathbb{R}^{n\times d}$ such that
\begin{equation}
\label{2.2} A_0 X_0=B_0.
\end{equation}
Here $A_0$ is the nonrandom true input matrix, $B_0$ is the true output
matrix, and $\tilde{A}$, $\tilde{B}$ are error matrices. The matrix
$X_0$ is the true value of the parameter.

We can rewrite the model \eqref{2.1}--\eqref{2.2} as a classical
functional errors-in-variables (EIV) model with vector regressor and
vector response \cite{ful}. Denote by $a_i^{\tr}$, $a_{0i}^{\tr}$,
$\tilde{a}_i^{\tr}$, $b_i^{\tr}$, $b_{0i}^{\tr}$, and $\tilde{b}_i^{\tr
}$ the rows of $A$, $A_0$, $\tilde{A}$, $B$, $B_0$, and $\tilde{B}$,
respectively, $i=1,\dots,m$. Then the model considered is equivalent to
the following EIV model:
\begin{equation*}
\label{2.3} a_i=a_{0i}+\tilde{a}_i, \qquad
b_i=b_{0i}+\tilde{b}_i, \qquad
b_{oi}=X_0^{\tr}a_{0i}, \quad i=1, \dots,m.
\end{equation*}
Based on observations $a_i$, $b_i$, $i=1,\dots,m$, we have to estimate
$X_0$. The vectors $a_{0i}$ are nonrandom and unknown, and the vectors
$\tilde{a}_i$, $\tilde{b}_i$ are random errors.

We state a global assumption of the paper.
\begin{enumerate}
\item\label{i}
The vectors $\tilde{z}_i$ with $\tilde{z}_i^{\tr}=[\tilde{a}_i^{\tr
},\tilde{b}_i^{\tr}]$, $i=1,2,\dots$, are i.i.d., with zero mean and
variance--covariance matrix
\begin{equation}
\label{2.4} S_{\tilde{z}}:=\cov(\tilde{z}_1)=
\sigma^2 \I_{n+d},
\end{equation}
where the factor of proportionality $\sigma^2$ is positive and unknown.
\end{enumerate}

The TLS problem consists in finding the values of disturbances $\Delta
\hat{A}$ and $\Delta\hat{B}$ minimizing the sum of squared corrections
\begin{equation}
\label{2.5} \min_{(X\in\mathbb{R}^{n\times d},\: \Delta A,\: \Delta B)}\bigl(\|\Delta A\|_F^2+\|
\Delta B\|_F^2\bigr)
\end{equation}
subject to the constraints
\begin{equation}
\label{2.6} (A-\Delta A)X=B-\Delta B.
\end{equation}

Here in \eqref{2.5}, for a matrix $C=(c_{ij})$, $\|C\|_F$ denotes the
Frobenius norm, $\|C\|_F^2=\sum_{i,j}c_{ij}^2$. Later on,
we will also use the operator norm $\|C\|=\sup_{x\neq
0}\tfrac{\|Cx\|}{\|x\|}$.
\subsection{TLS estimator and its consistency}
It may happen that, for some random realization, problem
\eqref{2.5}--\eqref{2.6} has no solution. In such a case, put $\hat
{X}_{\mathit{tls}}=\infty$. Now, we give a formal definition of the TLS estimator.
\begin{defin}
The TLS estimator $\hat{X}_{\mathit{tls}}$ of $X_0$ in the model \eqref
{2.1}--\eqref{2.2} is a~measurable mapping of the underlying
probability space into $\mathbb{R}^{n\times d}\cup\{\infty\}$, which
solves problem \eqref{2.5}--\eqref{2.6} if there exists a solution, and
$\hat{X}_{\mathit{tls}}=\infty$ otherwise.
\end{defin}

We need the following conditions for the consistency of $\hat{X}_{\mathit{tls}}$.
\addtocounter{enumi}{2}
\begin{enumerate}\addtocounter{enumi}{1}
\item\label{ii}
$\M\|\tilde{z}_1\|^4<\infty$, where $\tilde{z}_1$ satisfies condition
\eqref{i}.
\item\label{iii}
$\tfrac{1}{m}A_0^{\tr}A_0\to V_A$ as $m\to\infty$, where $V_A$ is a
nonsingular matrix.
\end{enumerate}

The next consistency result is contained in Theorem 4(a) of \cite{kukva}.
\begin{thm}\label{thm:2}
Assume condition \eqref{i} to \eqref{iii}. Then $\hat{X}_{\mathit{tls}}$ is
finite with probability tending to one, and $\hat{X}_{\mathit{tls}}$ tends to
$X_0$ in probability as $m\to\infty$.
\end{thm}
\subsection{The objective and estimating functions}
Denote\vspace{-6pt}
\begin{gather}
q(a,b;X)=\bigl(a^{\tr}X-b^{\tr}\bigr) \bigl(
\I_d+X^{\tr}X\bigr)^{-1}\bigl(X^{\tr}a-b
\bigr), \label
{2.7}
\\
Q(X)=\sum_{i=1}^m q(a_i,b_i;X),
\quad  X\in\mathbb{R}^{n\times d}. \label{2.8}
\end{gather}
The TLS estimator is known to minimize the objective function \eqref
{2.8}; see \cite{spr} or formula (24) in \cite{kukva}.
\begin{lemma}\label{l:3}
The TLS estimator $\hat{X}_{\mathit{tls}}$ is finite iff there exists an
unconstrained minimum of the function \eqref{2.8}, and then $\hat
{X}_{\mathit{tls}}$ is a minimum point of that function.
\end{lemma}

Introduce an estimating function related to the loss function \eqref{2.7}:
\begin{equation}
\label{2.9} s(a,b;X):=a\bigl(a^{\tr}X-b^{\tr}\bigr)-X\bigl(
\I_d+X^{\tr}X\bigr)^{-1}\bigl(X^{\tr}a-b
\bigr) \bigl(a^{\tr
}X-b^{\tr}\bigr).
\end{equation}
\renewcommand{\theenumi}{\alph{enumi}}
\renewcommand{\labelenumi}{\rm(\alph{enumi})}
\begin{cor}\label{cor:4}
\begin{enumerate}
\item\label{4.a}
Under conditions \eqref{i} to \eqref{iii}, with probability tending to
one $\hat{X}_{\mathit{tls}}$ is a solution to the equation\vspace{-3pt}
\begin{equation*}
\label{2.10} \sum_{i=1}^m
s(a_i,b_i;X)=0, \quad  X\in\mathbb{R}^{n\times d}.\vspace{-3pt}
\end{equation*}
\item\label{4.b}
Under assumption \eqref{i}, the function $s(a,b;X)$ is unbiased
estimating function, that is, for each $i\geq1$, $\M_{X_0}s(a_i,b_i;X_0)=0$.
\end{enumerate}
\end{cor}

Expression \eqref{2.9} as a function of $X$ is a mapping in $\mathbb
{R}^{n\times d}$. Its derivative $s_X^{\prime}$ is a linear operator in
this space.
\begin{lemma}\label{l:5}
Under condition \eqref{i}, for each $H\in\mathbb{R}^{n\times d}$ and
$i\geq1$, we have
\begin{equation}
\label{2.11} \M_{X_0}\bigl[s_X^{\prime}(a_i,b_i;X_0)
\cdot H\bigr]=a_{0i}a_{0i}^{\tr}H.
\end{equation}
\end{lemma}

Therefore, we can identify $\M_{X_0}s_X^{\prime}(a_i,b_i;X_0)$ with the
matrix $a_{0i}a_{0i}^{\tr}$.\vspace{-6pt}
\section{Main results}\label{s:3}
Introduce further assumptions to state the asymptotic normality of $\hat
{X}_{\mathit{tls}}$. We need a bit higher moments compared with conditions \eqref
{ii} and \eqref{iii} in order to use the Lyapunov CLT. Recall that
$\tilde{z}_i$ satisfies condition \eqref{i}.
\renewcommand{\theenumi}{\roman{enumi}}
\renewcommand{\labelenumi}{(\roman{enumi})}
\begin{enumerate} \addtocounter{enumi}{3}
\item\label{iv}
For some $\delta>0$, $\M\|\tilde{z}_1\|^{4+2\delta}<\infty$.
\item\label{v}
For $\delta$ from condition \eqref{iv},
\begin{equation*}
\label{3.1} \frac{1}{m^{1+\delta/2}}\sum_{i=1}^m\|a_{0i}\|^{2+\delta}
\to0 \quad  \text{ as } m\to\infty.\vspace{-6pt}
\end{equation*}
\item\label{vi}
$\frac{1}{m}\sum_{i=1}^m a_{0i}\to\mu_a$ as $m\to\infty$,
where $\mu_a\in\mathbb{R}^{n\times1}$.
\item\label{vii}
The distribution of $\tilde{z}_1$ is symmetric around the origin.
\end{enumerate}

Introduce a random element in the space of systems consisting of five matrices:
\begin{equation}
\label{3.2} W_i=\bigl(a_{0i}\tilde{a}_i^{\tr},a_{0i}
\tilde{b}_i^{\tr},\tilde{a}_i\tilde
{a}_i^{\tr}-\sigma^2\I_n,
\tilde{a}_i\tilde{b}_i^{\tr},
\tilde{b}_i\tilde {b}_i^{\tr}-
\sigma^2\I_d\bigr).
\end{equation}

Hereafter $\stackrel{\text{\rm d}}{\longrightarrow}$ stands for the
convergence in distribution.
\begin{lemma}\label{l:6}
Assume conditions \eqref{i} and \eqref{iii}--\eqref{vi}. Then
\begin{equation}
\label{3.3} \frac{1}{\sqrt{m}}\sum_{i=1}^m
W_i\stackrel{\text{\rm d}} {\longrightarrow} \varGamma= (
\varGamma_1,\dots,\varGamma_5) \quad \text{ as } m\to \infty,
\end{equation}
where $\varGamma$ is a Gaussian centered random element with matrix components.
\end{lemma}
\begin{lemma}\label{l:7}
In assumptions of Lemma~\ref{l:6}, replace condition \eqref{vi} with
condition \eqref{vii}. Then the convergence \eqref{3.3} still holds
with independent components $\varGamma_1,\dots,\varGamma_5$.
\end{lemma}

Now, we state the asymptotic normality of $\hat{X}_{\mathit{tls}}$.
\renewcommand{\theenumi}{\alph{enumi}}
\renewcommand{\labelenumi}{\rm(\alph{enumi})}
\begin{thm}\label{thm:8}
\begin{enumerate}
\item\label{8.a}
Assume conditions \eqref{i} and \eqref{iii}--\eqref{vi}. Then
\begin{equation}
\label{3.4} \sqrt{m}(\hat{X}_{\mathit{tls}}-X_0)
\stackrel{\text{\rm d}} {\longrightarrow}V_A^{-1}
\varGamma(X_0) \quad \text{as } m\to\infty,
\end{equation}
\begin{equation}
\label{3.5} \varGamma(X):=\varGamma_1 X - \varGamma_2
+ \varGamma_3 X - \varGamma_4 - X\bigl(\I_d
+ X^{\tr}X\bigr)^{-1}\bigl(X^{\tr}\varGamma_3
X - X^{\tr}\varGamma_4-\varGamma_4^{\tr}X
+ \varGamma_5\bigr),
\end{equation}
where $V_A$ satisfies condition \eqref{iii}, and $\varGamma_i$ satisfy
relation \eqref{3.3}.
\item\label{8.b}
In the assumption of part \eqref{8.a}, replace condition \eqref{vi}
with condition \eqref{vii}. Then the convergence \eqref{3.4} still
holds, and, moreover, the limit random matrix $X_{\infty
}:=V_A^{-1}\varGamma(X_0)$ has a nonsingular covariance structure, that
is, for each nonzero vector $u\in\mathbb{R}^{d\times1}$, $\cov
(X_{\infty}u)$ is a nonsingular matrix.
\end{enumerate}
\end{thm}
\begin{rem}
Conditions of Theorem~\ref{thm:8}\eqref{8.a} are similar to Gallo's
conditions \cite{gal} for the asymptotic normality in the univariate
case; see also, \cite{vaw}, pp. 240--243. Compared with Theorems 2.3
and 2.4 of \cite{pes}, stated for univariate case with mixing errors,
we need not the requirement for entries of the true input $A_0$ to be
totally bounded.
\end{rem}

In \cite{pes}, Section~2, we can find a discussion of importance of the
asymptotic normality result for $\hat{X}_{\mathit{tls}}$. It is claimed there
that the formula for the asymptotic covariance structure of $\hat
{X}_{\mathit{tls}}$ is computationally useless, but in case where the limit
distribution is nonsingular, we can use the block-bootstrap techniques
when constructing confidence intervals and testing hypotheses.

However, in the case of normal errors $\tilde{z}_i$, we can apply
Theorem~\ref{thm:8}\eqref{8.b} to construct the asymptotic confidence
ellipsoid, say, for $X_0 u$, $u\in\mathbb{R}^{d\times1}$, $u\neq0$.
Indeed, relations \cref{3.2,3.3,3.4,3.5} show that the nonsingular matrix
\[
S_u:=\cov\bigl(\mathbf{V}_A^{-1}
\varGamma(X_0)u\bigr)
\]
is a continuous function $S_u=S_u(X_0,\mathbf{V}_A,\sigma^2)$ of
unknown parameters $X_0$, $\mathbf{V}_A$, and $\sigma^2$. (It is
important here that now the components $\varGamma_j$ of $\varGamma$ are
independent, and the covariance structure of each $\varGamma_j$ depends on
$\sigma^2$ and $\mathbf{V}_A$, not on some other limit characteristics
of $A_0$; see Lemma 6.) Once we possess consistent estimators $\hat
{\mathbf{V}}_A$ and $\hat{\sigma}^2$ of $\mathbf{V}_A$ and $\sigma^2$,
the matrix
$\hat{S}_u:=S_u(\hat{X}_{\mathit{tls}}, \hat{\mathbf{V}}_A, \hat{\sigma}^2)$ is
a~consistent estimator for the covariance matrix $S_u$.

Hereafter, a bar means averaging for rows $i=1,\dots,m$, for example,
$\overline{ab^{\tr}}=m\frac{1}{m}\sum_{i=1}^{m}a_i b_i^{\tr}$.
\begin{lemma}\label{l:10}
Assume the conditions of Theorem~\ref{thm:2}. Define
\begin{equation}
\label{3.6} \hat{\sigma}^2 = \frac{1}{d}\mathrm{tr} \bigl[
\bigl(\overline{bb^{\tr}}- 2\hat{X}_{\mathit{tls}}^{\tr}
\overline{ab^{\tr}}+ \hat{X}_{\mathit{tls}}^{\tr}
\overline{aa^{\tr}}\hat{X}_{\mathit{tls}}\bigr) \bigl(\I_d+
\hat{X}_{\mathit{tls}}^{\tr}\hat{X}_{\mathit{tls}}\bigr)^{-1}
\bigr],
\end{equation}
\begin{equation*}
\label{3.7} \hat{V}_A = \overline{aa^{\tr}} - \hat{
\sigma}^2\I_n.
\end{equation*}
Then
\begin{equation}
\label{3.8} \hat{\sigma}^2\prob\sigma^2, \qquad
\hat{V}_A\prob V_A.
\end{equation}
\end{lemma}
\begin{rem}
Estimator \eqref{3.6} is a multivariate analogue of the maximum
likelihood estimator (1.53) in \cite{chn} in the functional scalar EIV model.
\end{rem}

Finally, for the case $\tilde{z}_1\sim N(0,\sigma^2\I_{n+d})$, based on
Lemma~\ref{l:10} and the relations
\[
\sqrt{m}(\hat{X}_{\mathit{tls}}-X_0 )\stackrel{\text{\rm
d}} {\longrightarrow}N(0,S_u), \quad S_u>0,\ \hat{S}_u\prob S_u,
\]
we can construct the asymptotic confidence ellipsoid for the vector
$X_0 u$ in a~standard way.
\begin{rem}
In a similar way, a confidence ellipsoid can be constructed for any
finite set of linear combinations of $X_0$ entries with fixed known
coefficients.
\end{rem}
%
\section{Conclusion}\label{s:4}
We extended the result of Gallo \cite{gal} and proved the asymptotic
normality of the TLS estimator in a multivariate model $A X\approx B$.
The normalized estimator converges in distribution to a random matrix
with quite complicated covariance structure. If the error distribution
is symmetric around the origin, then the latter covariance structure is
nonsingular. For the case of normal errors, this makes it possible to
construct the asymptotic confidence region for a vector $X_0 u$, $u\in
\mathbb{R}^{d\times1}$, where $X_0$ is the true value of $X$.

In future papers, we will extend the result for the elementwise
weighted TLS estimator \cite{kukva} in the model $A X\approx B$, where
some columns of the matrix $[A,B]$ may be observed without errors, and,
in addition, the error covariance matrix may differ from row to row.


\appendix
\section*{Appendix}
\section*{Proof of Corollary~\ref{cor:4}}

\renewcommand{\theequation}{4.\arabic{equation}}
(a) For any $n$ and $d$, the space $\mathbb{R}^{n\times d}$ is endowed
with natural inner product \linebreak$\langle A,B\rangle =\mathrm{tr}(AB^{\tr})$ and
the Frobenius norm. The matrix derivative $q_X^{\prime}$ of the
functional \eqref{2.7} is a linear functional on $\mathbb{R}^{n\times
d}$, which can be identified with certain matrix from $\mathbb
{R}^{n\times d}$ based on the inner product.

Using the rules of matrix calculus \cite{car}, we have for $H\in\mathbb
{R}^{n\times d}$:
\begin{align*}
\big\langle q_X^{\prime},H \big\rangle  &=a^{\tr}H\bigl(\I_d+X^{\tr}X \bigr)^{-1}\bigl(X^{\tr}a-b\bigr)\\
&\quad -\bigl(a^{\tr}X-b^{\tr}\bigr) \bigl(\I_d+X^{\tr}X \bigr)^{-1}\bigl(H^{\tr}X+X^{\tr}H\bigr) \bigl(\I_d+X^{\tr}X\bigr)^{-1}\bigl(X^{\tr}a-b \bigr)\\
&\quad +\bigl(a^{\tr}X-b^{\tr}\bigr) \bigl(\I_d+X^{\tr}X \bigr)^{-1}H^{\tr}a.
\end{align*}
Collecting similar terms, we obtain:
\begin{align*}
\label{4.1}
\frac{1}{2} \big\langle q_X^{\prime},H\big\rangle  &= \bigl(a^{\tr}X-b^{\tr}\bigr) \bigl(\I_d+X^{\tr}X \bigr)^{-1}H^{\tr}a\\
&\quad -\bigl(a^{\tr}X-b^{\tr}\bigr) \bigl(\I_d+X^{\tr}X\bigr)^{-1}H^{\tr}X\bigl(\I_d+X^{\tr}X\bigr)^{-1}\bigl(X^{\tr}a-b\bigr),
\end{align*}
and
\begin{align*}
\dfrac{1}{2} \big\langle q_X^{\prime},H \big\rangle  &= \mathrm{tr} \bigl[a\bigl(a^{\tr}X-b^{\tr}\bigr) \bigl(\I _d+X^{\tr}X\bigr)^{-1}H^{\tr} \bigr]\\
&\quad -\mathrm{tr} \bigl[X\bigl(\I_d+X^{\tr}X\bigr)^{-1}\bigl(X^{\tr}a-b\bigr) \bigl(a^{\tr}X-b^{\tr}\bigr) \bigl(\I _d+X^{\tr}X\bigr)^{-1}H^{\tr}\bigr].
\end{align*}

Using the inner product in $\mathbb{R}^{n\times d}$, we get $\tfrac
{1}{2}q_X^{\prime}=s(x)(\I_d+X^{\tr}X)^{-1}$, where $s(x)$ is the
left-hand side of \eqref{2.9}. In view of Theorem~\ref{thm:2} and Lemma
\ref{l:3}, this implies the statement of Corollary~\ref{cor:4}\eqref{4.a}.

(b) Now, we set
\begin{equation}
\label{4.2} a=a_0+\tilde{a}, \qquad b=b_0+\tilde{b}, \qquad
b_0=X^{\tr}a_0,
\end{equation}
where $a_0$ is a nonrandom vector, and, like in \eqref{2.4},
\begin{equation}
\label{4.3} \cov\lleft(\lleft[ %
\begin{array}{l}
\tilde{a}\\
\tilde{b}
\end{array} %
 \rright] \rright)= \sigma^2\I_{n+d}, \qquad \M\lleft[
\begin{array}{l}
\tilde{a}\\
\tilde{b}
\end{array} %
 \rright]=0.
\end{equation}
Then
\begin{equation}
\label{4.4} \M_X a\bigl(a^{\tr}X-b^{\tr}\bigr)=\M
a\bigl(\tilde{a}^{\tr}X-\tilde{b}^{\tr}\bigr)=\sigma
^2 X,
\end{equation}
\begin{equation}
\label{4.5} \M_X\bigl(X^{\tr}a-b\bigr)
\bigl(a^{\tr}X-b^{\tr}\bigr)=\M\bigl(X^{\tr}\tilde{a}-
\tilde {b}\bigr) \bigl(\tilde{a}^{\tr}X-\tilde{b}^{\tr}\bigr)=
\sigma^2\bigl(\I_d+X^{\tr}X\bigr).
\end{equation}
Therefore (see \eqref{2.9}),
\begin{equation*}
\M_X s(a,b;X)=\sigma^2 X - \sigma^2 X\bigl(
\I_d + X^{\tr}X\bigr)^{-1}\bigl(\I_d +
X^{\tr}X\bigr) = 0.
\end{equation*}
This implies the statement of Corollary~\ref{cor:4}\eqref{4.b}.
\section*{Proof of Lemma~\ref{l:5}}
The derivative $s_X^{\prime}$ of the function \eqref{2.9} is a linear
operator in $\mathbb{R}^{n\times d}$. For \linebreak$H\in\mathbb
{R}^{n\times d}$, we have:
\begin{equation}
\label{4.6} %
\begin{aligned}
s_X^{\prime}H&=aa^{\tr}H - H\bigl(\I_d + X^{\tr}X\bigr)^{-1}\bigl(X^{\tr}a - b\bigr) \bigl(a^{\tr}X - b^{\tr}\bigr)\\
&\quad  + X\bigl(\I_d + X^{\tr}X\bigr)^{-1}\bigl(H^{\tr}X + X^{\tr}H\bigr) \bigl(\I_d + X^{\tr}X\bigr)^{-1}\bigl(X^{\tr}a-b\bigr)\\
&\quad  \times\bigl(a^{\tr}X - b^{\tr}\bigr) - X\bigl(\I_d + X^{\tr}X\bigr)^{-1} \bigl(H^{\tr}a\bigl(a^{\tr}X - b^{\tr}\bigr) + \bigl(X^{\tr}a - b\bigr)a^{\tr}H\bigr).
\end{aligned} %
\end{equation}

As before, we set \eqref{4.2}, \eqref{4.3} and use relations \eqref
{4.4}, \eqref{4.5}, and the relation\break $\M aa^{\tr} = a_0a_0^{\tr} +
\sigma^2\I_n$. We obtain:
\begin{equation*}
\begin{aligned}
\M_X s_X^{\prime}H &= \bigl(a_0a_0^{\tr} + \sigma^2 \I_n\bigr)H - \sigma^2 H + \sigma^2 X\bigl(\I_d + X^{\tr}X\bigr)^{-1}\bigl(H^{\tr}X + X^{\tr}H\bigr)\\
& \quad - \sigma^2 X\bigl(\I_d + X^{\tr}X \bigr)^{-1}\bigl(H^{\tr}H + X^{\tr}H\bigr) = a_0 a_0^{\tr}H.
\end{aligned}
\end{equation*}
This implies \eqref{2.11}.
\section*{Proof of Lemma~\ref{l:6}}
The random elements $W_i$, $i\geq1$, in \eqref{3.2} are independent
and centered. We want to apply the Lyapunov CLT for the left-hand side
of \eqref{3.3}.

(a) All the second moments of $m^{-\frac{1}{2}}\sum_{i=1}^m W_i$ converge to finite limits. For example, for the first
component, we have
\begin{equation*}
\frac{1}{m}\sum_{i=1}^{m}\M
\bigl(\big\langle a_{0i} \tilde{a}_i^{\tr}, H_1\big\rangle
\bigr)^2 = \frac
{1}{m}\sum_{i=1}^{m}
\M\bigl(\text{tr } a_{0i} \tilde{a}_1^{\tr}
H_1^{\tr}\bigr)^2,
\end{equation*}
and this has a finite limit due to assumption \eqref{iii}. Here $H_1\in
\mathbb{R}^{n\times n}$, and we use the inner product introduced in the
proof of Corollary~\ref{cor:4}.

For the fifth component,
\begin{equation*}
\frac{1}{m}\sum_{i=1}^{m}\M\bigl(\big\langle
\tilde{b}_i\tilde{b}_i^{\tr} -
\sigma^2\I _d, H_2\big\rangle \bigr)^2 = \M
\bigl[\mathrm{tr}\bigl(\bigl(\tilde{b}_1\tilde{b}_1^{\tr}
- \sigma^2\I _d\bigr)H_2\bigr)
\bigr]^2 < \infty,
\end{equation*}
because the fourth moments of $\tilde{b}_i$ are finite. Here $H_2\in
\mathbb{R}^{d\times d}$.

For mixed moments of the first and fifth components, we have
\begin{align}
&\frac{1}{m}\sum_{i=1}^{m} \M \big\langle a_{0i}\tilde{a}_i^{\tr}, H_1\big\rangle \cdot \big\langle \tilde{b}_i\tilde{b}_i^{\tr} - \sigma^2\I_d, H_2\big\rangle \nonumber\\
&\quad = \M\Bigg\langle  \Biggl(\frac{1}{m}\sum_{i=1}^{m}a_{0i} \Biggr)\tilde{a}_1^{\tr}, H_1\Bigg\rangle \cdot\big\langle \tilde{b}_1\tilde{b}_1^{\tr} - \sigma^2\I_d, H_2\big\rangle ,\label{4.7} %
\end{align}
and this, due to condition \eqref{vi}, converges toward
\begin{equation*}
\M\big\langle \mu_a\tilde{a}_1^{\tr}, H_1\big\rangle
\cdot\big\langle \tilde{b}_1\tilde{b}_1^{\tr} -
\sigma^2\I_d, H_2\big\rangle .
\end{equation*}

Other second moments can be considered in a similar way.
\vspace{2mm}

(b) The Lyapunov condition holds for each component of \eqref{3.2}. Let
$\delta$ be the quantity from assumptions \eqref{iv}, \eqref{v}. Then
\begin{equation*}
\frac{1}{m^{1+\delta/2}}\sum_{i=1}^{m}
\M\big\|a_{0i}\tilde{a}_i^{\tr
}\big\|^{2+\delta}\leq
\frac{\M\|\tilde{a}_1\|^{2+\delta}}{m^{1+\delta/2}}\sum_{i=1}^{m}\|a_{0i}\|^{2+\delta}
\to0
\end{equation*}
as $m\to\infty$ by condition \eqref{v}. For the fifth component,
\begin{equation*}
\begin{aligned}
 \frac{1}{m^{1+\delta/2}}\sum_{i=1}^{m} \M\big\|\tilde{b}_i\tilde{b}_i^{\tr} - \sigma^2\I_d\big\|^{2+\delta} &= \frac{1}{m^{\delta/2}}\M\big\| \tilde{b}_1\tilde{b}_1^{\tr} - \sigma^2\I _d\big\|^{2+\delta}\\
&\leq\frac{\mathrm{const}}{m^{\delta/2}}\M\|\tilde{b}_1\|^{4+2\delta }\to0 \quad \text{ as } m\to\infty.
\end{aligned}
\end{equation*}
The latter expectation is finite by condition \eqref{iv}.

The Lyapunov condition for other components is considered similarly.
\vspace{2mm}

(c) Parts (a) and (b) of the present proof imply \eqref{3.3} by the
Lyapunov CLT.\vspace{-6pt}
\section*{Proof of Lemma~\ref{l:7}}
Under conditions \eqref{vii} and \eqref{i}, all the five components of
$W_i$, which is given in \eqref{3.2}, are uncorrelated (e.g., the
cross-correlation like \eqref{4.7} equals zero, and condition \eqref
{vi} is not needed). As in proof of Lemma~\ref{l:6}, the convergence
\eqref{3.3} still holds. The components $\varGamma_1,\dots,\varGamma_5$ of
$\varGamma$ are independent because the components of $W_i$ are uncorrelated.\vspace{-6pt}
\section*{Proof of Theorem~\ref{thm:8}(a)}
Our reasoning is typical for theory of generalized estimating
equations, with specific feature that a matrix parameter rather than
vector one is estimated.

By Corollary~\ref{cor:4}\eqref{4.a}, with probability tending to $1$ we have\vspace{-3pt}
\begin{equation}
\label{4.8} \sum_{i=1}^m
s(a_i,b_i;\hat{X}_{\mathit{tls}}) = 0.
\end{equation}

Now, we use Taylor's formula around $X_0$ with the remainder in the
Lagrange form; see \cite{car}, Theorem 5.6.2. Denote\vspace{-3pt}
\begin{equation*}
\hat{\Delta} = \sqrt{m}(\hat{X}_{\mathit{tls}} - X_0), \qquad
y_m = \sum_{i=1}^m
s(a_i,b_i;X_0), \qquad U_m = \sum
_{i=1}^m s_X^{\prime}(a_i,b_i;X_0).
\end{equation*}

Then \eqref{4.8} implies the relation\vspace{-3pt}
\begin{equation}
\label{4.9} %
\begin{aligned} & \biggl(\frac{1}{m}U_m
\biggr)\hat{\Delta} = - \frac{1}{\sqrt{m}}y_m + \mathit{rest}_1,
\\
&\|\mathit{rest}_1\|\leq\|\hat{\Delta}\|\cdot\|\hat{X}_{\mathit{tls}} -
X_0\|\cdot O_p(1). \end{aligned} %
\end{equation}
Here $O_p(1)$ is a factor of the form\vspace{-3pt}
\begin{equation}
\label{4.10} \frac{1}{m}\sum_{i=1}^m
\sup_{(\|X\|\leq\|X_0\|+1)}\big\|s_x^{\prime\prime
}(a_i,b_i;X)\big\|.
\end{equation}
Relation \eqref{4.9} holds with probability tending to $1$ because, due
to Theorem~\ref{thm:2}, $\hat{X}_{\mathit{tls}}\prob X_0$; expression \eqref
{4.10} is indeed $O_p(1)$ because the derivative $s_x^{\prime\prime}$
is quadratic in $a_i$, $b_i$ (cf. \eqref{4.6}), and the averaged second
moments of $[a_i^{\tr},b_i^{\tr}]$ are assumed to be bounded.

Now, $\|\mathit{rest}_1\|\leq\|\hat{\Delta}\|\cdot o_p(1)$. Next, by Lemma~\ref
{l:5} and condition \eqref{iii},
\begin{equation*}
\frac{1}{m}U_m = \frac{1}{m}\M U_m +
o_p(1) = V_A + o_p(1).
\end{equation*}
Therefore, \eqref{4.9} implies that\vspace{-6pt}
\begin{gather}
V_A\hat{\Delta} = -\frac{1}{\sqrt{m}}y_m +
\mathit{rest}_2,\label{4.11}
\\
\|\mathit{rest}_2\|\leq\|\hat{\Delta}\|\cdot o_p(1).\label{4.12}
\end{gather}

Now, we find the limit in distribution of $y_m/\sqrt{m}$. The summands
in $y_m$ have zero expectation due to Corollary~\ref{cor:4}\eqref{4.b}.
Moreover (see \eqref{2.9}),
\begin{equation*}
s(a_i,b_i;X_0)=(a_{0i} +
\tilde{a}_i) \bigl(\tilde{a}_i^{\tr}X_0
- \tilde {b}_i^{\tr}\bigr) - X_0\bigl(
\I_d + X_0^{\tr}X_0
\bigr)^{-1}\bigl(X_0^{\tr}\tilde{a}_i
- \tilde{b}_i\bigr) \bigl(\tilde{a}_i^{\tr}X_0
- \tilde{b}_i^{\tr}\bigr),\vspace{-3pt}
\end{equation*}
\begin{equation*}
\label{4.13} %
\begin{aligned}  s(a_i,b_i;X_0)&
= W_{i1}X_0 - W_{i2} + W_{i3}X_0
- W_{i4} - X_0\bigl(\I_d +
X_0^{\tr}X_0\bigr)^{-1}
\\
& \quad \times\bigl(X_0^{\tr}W_{i3}X_0 -
X_0^{\tr}W_{i4} - W_{i4}^{\tr}X_0
+ W_{i5}\bigr). \end{aligned} %
\end{equation*}
Here $W_{ij}$ are the components of \eqref{3.2}. By Lemma~\ref{l:6} we
have (see \eqref{3.5})
\begin{equation}
\label{4.14} \frac{1}{\sqrt{m}}y_m\stackrel{\text{\rm d}} {
\longrightarrow}\varGamma (X_0) \quad  \text{as } m\to\infty.
\end{equation}

Finally, relations \eqref{4.11}, \eqref{4.12}, \eqref{4.14} and the
nonsingularity of $V_A$ imply that $\hat{\Delta} = O_p(1)$, and by
Slutsky's lemma we get
\begin{equation}
\label{4.15} V_A\hat{\Delta}\stackrel{\text{\rm d}} {
\longrightarrow} \varGamma(X_0) \quad \text{as } m\to\infty.
\end{equation}
By condition \eqref{iii} the matrix $V_A$ is nonsingular. Thus, the
desired relation \eqref{3.4} follows from \eqref{4.15}.\vspace{-6pt}
\section*{Proof of Theorem~\ref{thm:8}(b)}
The convergence \eqref{3.4} is justified as before, but using Lemma~\ref
{l:7} instead of Lemma~\ref{l:6}. It suffices to show that $\cov(\varGamma
(X_0)u)$ is nonsingular for $u\in\mathbb{R}^{d\times1}$, $u\neq0$.

Now, the components $\varGamma_1,\dots,\varGamma_5$ are independent. Then\vspace{-3pt}
(see \eqref{3.5})
\begin{equation*}
\begin{aligned}
\cov\bigl(\varGamma(X_0)u\bigr)&\geq\cov(\varGamma_2 u) = \lim_{m\to\infty}\frac{1}{m}\sum_{i=1}^{m}\M\bigl(u^{\tr}\tilde{b}_i a_{0i}^{\tr}a_{0i}\tilde{b}_i^{\tr}u\bigr)\\
&= \mathrm{tr}V_A\cdot\M\big\|\tilde{b}_1^{\tr}u\big\|^2 = \sigma^2\rm{tr}V_A\cdot\|u\|^2>0.
\end{aligned} %
\end{equation*}
\section*{Proof of Lemma~\ref{l:10}}
By condition \eqref{i} we have\vspace{-6pt}
\begin{equation*}
\begin{aligned} \M a_i a_i^{\tr} &=
a_{0i}a_{0i}^{\tr} + \sigma^2
\I_n, \qquad \M a_i b_i^{\tr} =
a_{i0}a_{i0}^{\tr}X_0,
\\
\M b_i b_i^{\tr} &= X_0^{\tr}a_{0i}a_{0i}^{\tr}X_0
+ \sigma^2\I_d, \end{aligned} %
\end{equation*}
\begin{equation}
\label{4.16} \M b_i b_i^{\tr} -
2X_0^{\tr}\M a_i b_i^{\tr}
+ X_0^{\tr}\bigl(\M a_i a_i^{\tr}
\bigr)X_0 = \sigma^2\bigl(\I_d +
X_0^{\tr}X_0\bigr).
\end{equation}
Equality \eqref{4.16} implies the first relation in \eqref{3.8} because
$\hat{X}_{\mathit{tls}}\prob X_0$ and $\overline{aa^{\tr}} - \M\overline{aa^{\tr
}}\prob0$, $\overline{ab^{\tr}} - \M\overline{ab^{\tr}}\prob0$,
$\overline{bb^{\tr}} - \M\overline{bb^{\tr}}\prob0$,

Finally,
\begin{align*}
&\hat{V}_A = \M\overline{aa^{\tr}} + o_p(1) - \hat{\sigma}^2\I_n = \overline{a_0a_0^{\tr}} + (\sigma^2 - \hat{\sigma}^2)\I_n + o_p(1),\\
&\quad \hat{V}_A\prob\displaystyle\lim_{m\to\infty}\overline{a_0a_0^{\tr}} = V_A.
\end{align*}

%

\end{document}